\documentclass[a4paper,11pt,draft]{article}
\usepackage[english,francais]{babel}
\usepackage{amsfonts,amsmath}
\newtheorem{thm}{Th{\'e}or{\`e}me}

\newtheorem{conj}[thm]{Conjecture}
\newtheorem{cor}[thm]{Corollaire}
\newtheorem{prop}[thm]{Proposition}
\newtheorem{rem}[thm]{Remarque}

\newcommand{\Z}{\mathbb{Z}}
\newcommand{\prodd}{\mathop{\prod}\limits}
\newcommand{\summ}{ \mathop{\sum}\limits}

\begin{document}
\title{D{\'EVELOPPEMENTS LIMIT\'ES ET LA TRANSFORM\'EE INVERSE}}
\author{Roland BACHER}
\date{}
\maketitle

{\def\thefootnote{\relax}
} 

\selectlanguage{english}
\begin{abstract}
 The generating series associated to a certain sequence of limited Taylor expansions coincides with the continuously iterated Invert Transform and displays
thus a nice functional equality.
\end{abstract}
\selectlanguage{francais}
\begin{abstract} 
La s\'erie g\'en\'eratrice d'une certaine suite de d\'eveloppements
limit\'es co\"{\i}ncide avec l'it\'eration continue de la transform\'ee
inverse et v\'erifie donc une jolie \'equation fonctionnelle.
\end{abstract}

\bigskip
\section{Le r\'esultat principal}

Notons $\lfloor f(x)\rfloor_k=\summ_{j=0}^k \varphi_j x^j$ le d{\'e}veloppement
limit{\'e} d'ordre $k$ d'une s{\'e}rie formelle $f(x)=\summ_{j=0}^\infty
\varphi_jx^j$. Associons \`a une s\'erie formelle $s(x)=1+\sum_{j=1}^\infty s_jx^j$ de
coefficient constant
$s_0=1$ la suite de polyn{\^o}mes
$$P_0(x)=1,\ P_1(x)=1+s_1x,\ P_2(x)=1+2s_1x+(s_1^2+s_2)x^2,\dots$$
d{\'e}finie de fa\c con r{\'e}cursive par 
$$P_0(x)=s_0=1 \mbox{~~et~~} P_k(x)=\lfloor 
P_{k-1}(x)\ s(x))\rfloor_k$$
et posons $Q_n(x)=x^nP_n(1/x)$. Le r\'esultat principal de ce travail est
une \'equation fonctionnelle pour la s\'erie g\'en\'eratrice
$\sum_{n=0}^\infty Q_n(x)t^n$ de la suite $Q_0(x),Q_1(x),\dots$:

\begin{thm}\label{thmA} On a 
$$\sum_{n=0}^\infty Q_n(x)t^n=\frac{\sum_{n=0}^\infty Q_n(0)t^n}{1-tx\sum_{n=0}^\infty Q_n(0)t^n}\ .$$
\end{thm}

La prochaine section de ce papier donne une interpr\'etation 
de cette \'equa\-tion fonctionnelle en termes de la transform\'ee
inverse d'une suite.

La derni\`ere section contient l'\'ebauche d'une preuve,
bas\'ee essentiellement sur une identit\'e de Stam \cite{S} qui
g\'en\'eralise un r\'esultat de Hurwitz \cite{Hu}.

\section{La transform\'ee inverse}
On associe {\`a} la suite $a=(a_0,a_1,a_2,\dots)$ ({\`a} valeurs dans un 
anneau commutatif quelconque) sa {\em transform{\'e}e inverse}
(Invert transform) $I(a)=b=(b_0,b_1,\dots)$ d{\'e}finie formellement par 
l'{\'e}galit{\'e}
$$\left(1+t\summ_{n=0}^\infty a_nt^n\right)\left(1-t\summ_{n=0}^\infty 
b_nt^n\right)=1\ .$$
Posons $I^0(a)=a,\ I^1(a)=b$ et $I^{k+1}(a)=I(I^k(a))$. Comme 
$a\longmapsto I(a)$ est bijective sur l'ensemble des suites, 
on peut donner un sens 
{\`a} $I^k(a)$ pour tout $k\in\Z$ en posant
$I^k(I^{-k}(a))=a$. Un petit calcul et une r{\'e}currence
montrent facilement le r{\'e}sultat suivant.

\begin{prop}\label{invtransform} Soit $k\in \Z$ un entier et soit
$a=(a_0,a_1,\dots)$ une suite. La s\'erie g\'en\'eratrice 
$\summ_{n=0}^\infty b_n t^n$ de la suite
$I^k(a)=b=(b_0,b_1,b_2,\dots)$ v{\'e}rifie alors l'identit\'e
$$\summ_{n=0}^\infty b_n t^n=\frac{\summ_{n=0}^\infty a_n t^n}
{1+kt\summ_{n=0}^\infty a_n t^n}\ .$$
\end{prop}
Cette proposition permet d'interpoler les it\'er\'ees de la 
transform{\'e}e inverse
$$I^x(a)=(I_0(x)=a_0,\ I_1(x)=a_1-a_0^2x,\ I_2(x)=a_2-2a_0a_1x+a_0^3x^2
,\dots$$ 
d'une suite $a=(a_0,a_1,a_2,\dots)$ en des valeurs
arbitraires . La notation $I^x(a)$ est justifi{\'e}e par l'identit{\'e}
$I^x(I^y(a))=I^{x+y}(a)$ et le th\'eor\`eme \ref{thmA} peut donc aussi 
s'\'enoncer sous la forme
$$I^x(Q(0))=Q(-x)$$
o\`u $Q(x)$ d\'esigne la suite $Q_0(x),Q_1(x),\dots$.

\begin{rem}\em{Un ph{\'e}nom{\`e}ne similaire d'interpolation continue se produit
{\'e}galement pour la composition it{\'e}r{\'e}e $f^{\circ k}=f\circ f\circ
\dots\circ f$ d'une s{\'e}rie formelle $f(t)=t+\summ_{i=2}^\infty a_i t^i$
dont le d{\'e}veloppement {\`a} l'ordre $1$ est l'identit{\'e} (ceci se 
g{\'e}n{\'e}ralise d'ailleurs facilement {\`a} un d$-$uplet de s{\'e}ries formelles 
$F(t_1,\dots,t_d)=(f_1(t_1,\dots,t_d),\dots,f_d(t_1,\dots,t_d))$
v{\'e}rifiant $f_i=t_i+$ termes d'ordre $>1$). Il existe alors une suite
$$\begin{array}{l}
C_1(x)=1,\ C_2(x)=a_2x,\ C_3(x)=(a_2^2(x-1)+a_3)x,\\
\quad C_4(x)=(((2x-3)a_2^3+5a_2a_3)(x-1)+2a_4)x/2,\dots\end{array}$$
avec $C_n(x)$ polynomial de degr{\'e} $\leq n-1$ en $x$ telle qu'on ait
$f^{\circ x}(t)=\summ_{i=1}^\infty C_i(x)t^i$.

Pour le prouver on peut consid{\'e}rer la diff{\'e}rence finie
$$C_n(k+1)-C_n(k)=\hbox{coefficient de }t^n\hbox{ dans }
\summ_{i=2}^{\infty}a_i\left(\summ_{j=1}^{\infty}C_j(k)t^j
\right)^i$$
qui est polynomiale de degr{\'e} au plus $n-2$ en $k$ par r{\'e}currence sur 
$n$. On peut {\'e}galement le d{\'e}duire de l'existence d'un isomorphisme de 
groupe entre ces s{\'e}ries (avec pour produit la composition $f\circ g$) 
et un certain groupe de matrices (infinies) 
triangulaires sup{\'e}rieures unipotentes, cf. Theorem 1.7a dans \cite{H}.}
\end{rem}

\begin{rem}\label{toeplitz} \em{Signalons encore la propri{\'e}t{\'e} suivante de
la transform{\'e}e inverse: {\'e}tendons la suite $b=I(a)$ en posant $b_{-1}=-1$
et $b_{-k}=0$ pour $k\geq 2$. Alors $a_n=\det\left((b_{i-j})_{0\leq i,j<n})\right)$ (la preuve est donn{\'e}e par l'isomorphisme entre l'anneau des s{\'e}ries
formelles et l'alg{\`e}bre des matrices de Toeplitz triangulaires).}
\end{rem}

Notons $(1^{\mu_1}\cdot 2^{\mu_2}\cdots m^{\mu_m})$
la partition de l'entier naturel
$m=\summ_{j=1}^m j\mu_j$ ayant $\mu_j$ parts de longueur $j$ et
d\'esignons par ${\cal P}_m$ l'ensemble fini de toutes
les partitions de $m$. 
Rappelons {\'e}galement la d{\'e}finition des coefficients multinomiaux
$${n\choose \nu}={n\choose \nu_1,\nu_2,\dots}=\frac{n!}{\Bigl(\prodd_j \nu_j!\Bigr)(n-\summ_j \nu_j)!}$$
pour $\nu=(\nu_1,\nu_2,\dots)$ une suite finie d'entiers naturels 
de somme $\summ \nu_i\leq n$.
En posant $a^\nu=\prodd_{j=1} a_j^{\nu_j}$ et en appliquant le th{\'e}or{\`e}me 
binomial {\`a} $\summ_{n=0}^\infty I_n(x)t^n=\summ_{n=0}^\infty
\Bigl(\summ_{j=0}^\infty a_j t^j\Bigr)^{n+1}(-xt)^n$ on montre:

\begin{prop}\label{formuleI}
Pour $a=(a_0=1,a_1,a_2,\dots)$ le polyn{\^o}me
$I_n(x)$ est donn{\'e} par la formule
$$I_n(x)=\summ_{k=0}^n (-x)^k
\summ_{(\nu_1,\nu_2,\dots),\ (1^{\nu_1}\cdot2^{\nu_2}\cdots)
\in {\cal P}_{n-k}}{1+k\choose \nu}a^\nu\ .$$
\end{prop}

Une formule dans le cas g{\'e}n{\'e}ral ($a_0\not=1$) se d{\'e}duit 
des propri{\'e}t{\'e}s d'homog{\'e}n{\'e}it{\'e} du polyn{\^o}me $I_n(x)$.

La $n-$i{\`e}me \em{matrice de Hankel} $H(n)$ d'une suite 
$s=(s_0,s_1,s_2,\dots)$ est la matrice avec coefficients 
$h_{i,j}=s_{i+j},\ 0\leq i,j<n$.
La \em{transform{\'e}e de Hankel} de $s$ est alors d{\'e}finie
comme {\'e}tant la suite
$$\det(H(1)),\det(H(2)),\det(H(3)),\dots\ $$
des d{\'e}terminants des matrices de Hankel d'ordre $1,2,\dots$ 
associ{\'e}es {\`a} $s$.

Layman \cite{L} a montr{\'e} le r{\'e}sultat suivant:

\begin{thm}\label{layman} (Layman) Deux suites $a$ et $b=I(a)$ 
reli{\'e}es par la transform{\'e}e inverse ont m{\^e}me transform{\'e}e de Hankel.
\end{thm}

Comme les polyn{\^o}mes $I_n(x)$ interpolent les transform{\'e}es 
inverses it{\'e}r{\'e}es, on a:

\begin{cor}\label{laymancontinue} La transform{\'e}e de Hankel de la suite
$I^x(a)=(I_0(x),I_1(x)$, $I_2(x),\dots)$
ne d{\'e}pend pas de $x$.
\end{cor}

Pour $k>0$ un entier, d{\'e}finissons la $k-$i{\`e}me 
transform{\'e}e de Hankel de  $s=(s_0$, $s_1,\dots)$ comme
la suite $\big(\det(H_k(n))\big)_{n=1,2,\dots}$ pour
$H_k(n)=(s_{i+j+k})_{0\leq i,j<n}$. 

\begin{conj}\label{laymangen} (i) La suite
$$\Big(\det\big((I_{i+j+k}(x))_{0\leq i,j<n}\big)\Big)_{n=1,2,3,\dots}$$
de la $k-$i{\`e}me transform{\'e}e de Hankel de $I^x(a)=(I_0(x),I_1(x),\dots)$
ne contient que des polyn{\^o}mes de degr{\'e} $\leq k$ en $x$.

\ \ (ii) Le d\'eterminant 
$$\det((Q_{i+j}(x))_{0\leq i,j<n})$$
pour $Q_0,Q_1,\dots$ associ\'es \`a $s(x)=1+s_1x+s_2x^2+\dots$ comme dans 
la section pr\'ec\'edente, ne d\'epend pas de $s_1$.
\end{conj}

\section{Id\'ee de la preuve du th\'eor\`eme \ref{thmA}}

L'ingr\'edient principal de la preuve du th\'eor\`eme \ref{thmA}
est le r\'esultat suivant qui exprime les coefficients des polyn\^omes
$P_i(x)$ (ou des polyn\^omes $Q_i(x)$) en fonctions des coefficients
$s_0=1,s_1,s_2,\dots$ de la s\'erie formelle $s(x)=\sum s_ix^i$ 
de d\'epart. 

\begin{thm}\label{formuleP} Pour $s(x)=\sum s_ix^i=1+s_1x+s_2x^2+\dots$
on a
$$P_n(x)=\summ_{k=0}^n x^k\summ_{\nu=(1^{\nu_1}\cdot 2^{\nu_2}\cdots)
\in{\cal P}_k}\frac{n+1-k}{n+1-\summ_j \nu_j}{n\choose \nu}s^\nu\ .$$
\end{thm}

\noindent{\bf Preuve du th\'eor\`eme \ref{thmA}.} En appliquant 
le th\'eor\`eme \ref{formuleP} aux deux expressions
$$\left(\summ_j Q_j(x)t^j\right)\left(
1-tx\summ_j Q_j(0)t^j\right)\quad \hbox{ et }\quad\summ_j Q_j(0)t^j$$
et on comparant  les coefficients de $x^nt^ka^\nu$, on
est ramen\'e \`a un cas particulier de l'identit\'e (10) dans le papier
\cite{S} de Stam.\hfill QED

\vspace{15pt}
\noindent{\bf Id{\'e}e de la preuve du th{\'e}or{\`e}me \ref{formuleP}.} {\`A} une
partition $\nu=(1^{\nu_1}\cdot 2^{\nu_2}\cdots)$ on associe le
polyn{\^o}me 
$$\begin{array}{ccl}
R_{\nu}(x)&=&\frac{x+1-\summ_j j\nu_j}{x+1-\summ_j \nu_j}{x\choose \nu}\\
&=&(x+1-\summ j\nu_j){x\choose (\summ \nu_j)-1}
\frac{((\summ \nu_j)-1)!}{\prodd \nu_j!}\end{array}\ .$$
Un calcul montre que les polyn{\^o}mes $R_\nu(x)$ satisfont
l'identit{\'e}
$$R_{\nu}(x+1)-R_{\nu}(x)=\summ_{j,\ \nu_j>0}R_{(1^{\nu_1}\cdots
(j-1)^{\nu_{j-1}}\cdot j^{\nu_j-1}\cdot (j+1)^{\nu_{j+1}}\cdots)}(x)\ .$$
On a ensuite
$$P_{k+1}(x)-P_k(x)=\lfloor P_k(x)\left(\summ_{j=1}a_j x^j\right)\rfloor_{k+1}$$
par d{\'e}finition et les polyn{\^o}mes
$$\tilde P_k(x)=1+\summ_{j=1}^k \summ_{(1^{\nu_1}\cdot 2^{\nu_2}\cdots)\in
{\cal P}_j}\left(\prodd_{s=1} a_s^{\nu_s}\right)\ R_{\nu}(k)\ x^j$$
satisfont la m{\^e}me {\'e}quation. Par r{\'e}currence sur $k$, on a donc
$\tilde P_k(x)-P_k(x)=c_kx^k$ avec $c_k$ une constante. Or $c_k$ doit {\^e}tre 
nul car tous les polyn{\^o}mes $R_\nu(x),\ \nu\in {\cal P}_k$ contribuant
au coefficient dominant $x^k$ de $\tilde P_k(x)$ 
ont une racine commune en $k-1$.\hfill QED

\vspace{20pt}
\emph{Je remercie R. Chapman qui m'a communiqu\'e une preuve de cas
particulier de l'identit\'e de Stam n\'ecessaire \`a la preuve du
th\'eor\`eme \ref{thmA} et I. Gessel qui m'a signal\'e
le travail de Stam.}

\vfill
\begin{flushleft}
Roland BACHER\\
INSTITUT FOURIER\\
Laboratoire de Math{\'e}matiques\\
UMR5582 (UJF-CNRS)\\
BP 74\\
38402 St MARTIN D'H{\`E}RES Cedex (France)\\[10pt]
Roland.Bacher@ujf-grenoble.fr     
\end{flushleft}

\end{document}